\documentclass{article}

\usepackage{amssymb}

\newtheorem{theorem}{Theorem}[section]
\newtheorem{lem}[theorem]{Lemma}
\newtheorem{prop}[theorem]{Proposition}

\newtheorem{dfn}[theorem]{Definition}
\newtheorem{prob}[theorem]{Problem}

\newcommand\A{\mathcal{A}}
\newcommand\B{\mathcal{B}}

\newcommand\U{\mathcal{U}}

\renewcommand\O{\Omega}

\newcommand\bb[1]{\mathbb{#1}}

\newcommand\mod{\mathrm{mod}\,}
\newcommand\qed{\hspace{\stretch{1}}$\Box$}

\begin{document}

\title{A pair of non-homeomorphic product measures on the Cantor set}
\author{Tim D. Austin}
\date{}

%\thanks This work was carried out during a summer research studentship
%funded by Trinity College, Cambridge over the 2005 summer
%vacation.}

\maketitle

\begin{abstract}
For $r \in [0,1]$ let $\mu_r$ be the Bernoulli measure on the
Cantor set given as the infinite power of the measure on $\{0,1\}$
with weights $r$ and $1-r$.  For $r,s \in [0,1]$ it is known that
the measure $\mu_r$ is continuously reducible to $\mu_s$ (that is,
there is a continuous map sending $\mu_r$ to $\mu_s$) if and only
if $s$ can be written as a certain kind of polynomial in $r$; in
this case $s$ is said to be binomially reducible to $r$.  In this
paper we answer in the negative the following question posed by
Mauldin:

\emph{Is it true that the product measures $\mu_r$ and $\mu_s$ are
homeomorphic if and only if each is a continuous image of the
other, or, equivalently, each of the numbers $r$ and $s$ is
binomially reducible to the other?}
\end{abstract}

\section{Introduction}

Two Borel measures $\mu$ and $\nu$ on a topological space $\O$ are
said to be \textbf{homeomorphic} if there is some
autohomeomorphism $h$ of the underlying space $\O$ such that $\mu$
is mapped to $\nu$ under $h$: $\nu = \mu \circ h^{-1}$. This means
that $\nu(A) = \mu(h^{-1}(A))$ for any Borel $A \subseteq \O$.

Characterizations of when measures are equivalent under
homeomorphisms have been given for a variety of special
topological spaces: for the $n$-dimensional unit cube by Oxtoby
and Ulam \cite{OxUl}, for the irrationals in the unit interval
(that is, the Baire space $\bb{N}^\bb{N}$) by Oxtoby \cite{Ox} and
for the Hilbert cube by Oxtoby and Prasad \cite{OxPras}. In this
paper we consider the Cantor set $\O = \{0,1\}^{\bb{N}}$ with its
Cartesian product topology, and restrict attention to those
probability measures $\mu$ on $\O$ that are given by an infinite
power of a probability measure $\lambda$ on $\{0,1\}$. These are
the measures arising in a description of a sequence of independent
tosses of a biased coin. For $r \in [0,1]$ we will write $\mu_r$
for the infinite power of the measure $\lambda_r$ with
$\lambda_r\{1\} = r$, $\lambda_r\{0\} = 1-r$.

If we ask first when one measure $\mu_s$ is the continuous image
of another, say $\mu_r$, in the sense that there is a continuous
self-map $f$ of $\O$ with $\mu_s = \mu_r \circ f^{-1}$, we find
that it is equivalent to an algebraic condition on $r$ and $s$:
that $s$ be \textbf{binomially reducible} to $r$.  This is defined
later. Given the definition, it is now easy to check the following
proposition (or see \cite{Mauld}):

\begin{prop}
For $r,s \in [0,1]$, $\mu_s$ is the continuous image of $\mu_r$ if
and only if $s$ is binomially reducible to $r$. \qed
\end{prop}

An amusing consequence of this is that binomial equivalence is a
genuine equivalence relation.

Homeomorphisms of power measures on the Cantor set have been
studied by Navarro-Berm\'udez \cite{Nav}, where it was proved that
for $r$ rational or transcendental $\mu_r$ is homeomorphic to
$\mu_s$ only if $s$ is $r$ or $1-r$, and where binomial
equivalence of $r$ and $s$ as a necessary condition for $\mu_r$ to
be homeomorphic to $\mu_s$ was obtained. Huang \cite{Hu} showed
that for $r$ an algebraic integer of degree 2 we are still
restricted to the trivial cases $s = r,1-r$, but constructed
examples of non-trivial $r$ and $s$ of larger degree which are
binomially equivalent. In \cite{NavOx}, Navarro-Berm\'udez and
Oxtoby showed that one of these examples studied by Huang does in
fact give a pair of homemomorphic measures.  More recently,
Dougherty, Mauldin and Yingst \cite{MauldDoughYingst} have found a
proof that yields many more examples in a systematic way.

The following problem appears in \cite{Mauld} as Problem 1065:

\begin{prob}
Is it true that the product measures $\mu_r$ and $\mu_s$ are
homeomorphic if and only if each is a continuous image of the
other, or, equivalently, the numbers $r$ and $s$ are binomially
equivalent.
\end{prob}

In this paper we will construct an example to show that the answer
to this question is No.

In the rest of the Introduction we will explain binomial
reducibility and introduce some related concepts that will be
needed later on.

First of all we fix some notation and terminology.  Let $e$ be a
sequence of 0s and 1s indexed by some finite subset $S$ of
$\bb{N}$.  We write $\langle e \rangle$ for the set of all
sequences in $\O$ whose $i^{\rm{th}}$ term agrees with that of $e$
for $i \in S$. The product topology on $\O$ has a base consisting
of clopen sets of the form $\langle e \rangle$. We will refer to
such sets as \textbf{cylinders}, and will say that $\langle e
\rangle$ has \textbf{length} $n$ if $n = \#S$ and that $\langle e
\rangle $ \textbf{depends on} $S$.  For a fixed $i$ we also write
$A_i$ for the cylinder given by specifying only that the
$i^{\rm{th}}$ term of a sequence be 1.

We now make some definitions.

\begin{dfn}
A polynomial $P(X)$ with integer coefficients will be called a
\textbf{partition polynomial} if it can be written in the form
\[a_nX^n + a_{n-1}X^{n-1}(1-X) + \cdots + a_0(1-X)^n\]
for some $n \geq 0$ and some integers $a_0,a_1,\ldots,a_n$
satisfying $0 \leq a_i \leq {n \choose i}$ for all $i \leq n$. We
will call this a \textbf{partition form} for $P$. The
\textbf{depth} of $P$ is the least $n$ for which such a partition
form exists.
\end{dfn}

Although we will not use their result, we note that Dougherty,
Mauldin and Yingst have recently given a simple characterization
of those integer polynomials that are partition polynomials; see
\cite{MauldDoughYingst}.

We are now in a position to define binomial reducibility:

\begin{dfn}
Given $r,s \in [0,1]$, we say that $s$ is \textbf{binomially
reducible} to $r$ if there is a partition polynomial $P$ with $s =
P(r)$.  We say that $r$ and $s$ are \textbf{binomially equivalent}
if each is binomially reducible to the other.
\end{dfn}

The depth of a partition polynomial is clearly no less than its
degree, but it can be strictly more: the polynomial $3X(1-X)$ is a
partition polynomial of degree 2, but to put it into partition
form we must write it as $3X(1-X)^2 + 3X^2(1-X)$ and so its depth
is 3.

We can see that if a partition polynomial has a partition form
with a given value of $n$ then it also has one with any larger
value of $n$ (by multiplying through by $X + (1-X)$).

Partition polynomials relate to power measures on the Cantor set
in the following way. If $\langle e \rangle$ is a cylinder where
$e$ has $a$ terms equal to 1 and $b$ terms equal to 0 then
(directly from the definition of product measure) we have
$\mu_r(\langle e \rangle) = r^a(1-r)^b$. The clopen subsets of
$\O$ are precisely those that can be written as a finite disjoint
union of cylinders (informally, those that depend on only finitely
many coordinates), and so, by breaking these cylinders up into
smaller cylinders as necessary, any such set $A$ can actually be
written as a finite disjoint union of cylinders all depending on
the same finite set of coordinates, say $S$ with $\# S = n$ . It
follows that $\mu_r(A)$ is equal to a finite sum of terms of the
form $r^{n-a}(1-r)^a$, as this is the measure of any cylinder
$\langle e \rangle$ where $e$ depends on $S$ and has $a$ 0s and
$n-a$ 1s. Since there are at most ${n \choose a}$ such cylinders,
the multiplicity of this term in the sum must be at most ${n
\choose a}$.  Thus we see that the values of the form $\mu_r(A)$
for $A$ clopen are precisely the values $P(r)$ for $P$ a partition
polynomial.

\begin{dfn}
A partition polynomial $P$ is said to \textbf{represent} a clopen
subset $A$ of $\O$ if $\mu_r(A) = P(r)$ for all $r \in [0,1]$.
\end{dfn}

A partition polynomial $P$ represents $A$ if and only there is a
partition form for $P$, say
\[P(X) = a_nX^n + a_{n-1}X^{n-1}(1-X) + \cdots + a_0(1-X)^n,\]
such that we can write $A$ as the disjoint union of cylinders
$U_1,U_2,\ldots,U_N$, all of length $n$, where $N = a_n + a_{n-1}
+ \cdots + a_0$ and for each $i \leq n$, precisely $a_i$ of the
sets $U_j$ have $i$ coordinates specified to be 1 and $n-i$
specified to be 0.  This is clearly sufficient, for if this latter
condition holds then
\[\mu_r(A) = \sum_{j \leq N}\mu_r(U_j),\]
and for each $i \leq n$ precisely $a_i$ of the terms of this sum
are equal to $r^i(1-r)^{n-i}$.  To see that the condition is
necessary, we need only choose $n$ to be at least the depth of $P$
and also so big that $A$ can be written as a disjoint union of
cylinders of length $n$, and then observe that both $P(r)$ and
$\mu_r(A)$ can be written as sums of terms of the form
$r^i(1-r)^{n-i}$; since these are linearly independent as
functions of $r \in [0,1]$ the coefficients must agree and the
result follows.

The same reasoning shows that any polynomial $P$ with integer
coefficients is a partition polynomial if and only if there is a
clopen subset $A$ of $\O$ such that $P(r) = \mu_r(A)$ for all $r
\in [0,1]$.

Later we will need the following lemma.

\begin{lem}\label{lem:polycomp}
If $P$ and $Q$ are partition polynomials then so is $P \circ Q$
\end{lem}

\noindent\textbf{Proof}\hspace{5pt} Choose a clopen set $A$
represented by $Q$ and a clopen set $B$ represented by $P$. By the
above comments it will suffice to find a clopen set $C$ such that
$\mu_r(C) = P(Q(r))$ for all $r \in [0,1]$.

Suppose we can write $B$ as a disjoint union of cylinders of a
fixed length $n$.  Then there is a subset $\B$ of $\{0,1\}^n$ such
that
\[B = \bigcup_{e \in \B}\langle e \rangle.\]
Similarly, suppose $A$ depends only on the first $m$ coordinates,
and write
\[A = \bigcup_{\eta \in \A}\langle \eta \rangle\]
for some $\A \subseteq \{0,1\}^m$.  Now define
$A_1,A_2,\ldots,A_n$ to be independent copies of $A$, in the
following sense:
\begin{eqnarray*}
A_1 = A &=& \{x \in \O:\ (x_1,x_2,\ldots,x_m) \in \A\}\\ A_2 &=&
\{x \in \O:\ (x_{m+1},x_{m+2},\ldots,x_{2m}) \in \A\}\\ &\vdots&\\
A_n &=& \{x \in \O:\ (x_{nm-m+1},x_{nm-m+2},\ldots,x_{nm}) \in
\A\}
\end{eqnarray*}
For each $e = (e_1,e_2,\ldots,e_n) \in \B$ let $C_e$ be the set
$C_{e,1} \cap C_{e,2} \cap \cdots \cap C_{e,n}$ where
\begin{eqnarray*}
C_{e,i} = \left\{\begin{array}{ll} A_i &\hbox{if $e$ has
$i^{\rm{th}}$ term equal to 1}\\ \O \setminus A_i &\hbox{if $e$
has $i^{\rm{th}}$ term equal to 0} \end{array}\right.
\end{eqnarray*}
Finally, let $C$ be $\bigcup_{e \in \B}C_e$.  It is easy to see
that the sets $C_e$ are disjoint (since for any two of them there
will be some $i \leq n$ such that one is contained in $A_i$ and
the other in $\O\setminus A_i$). Also,
\[\mu_r(C) = \sum_{e \in \B}\mu_r(C_e),\]
and if $e$ has $i$ 1s and $n-i$ 0s then $\mu_r(C_e) = \mu_r(A)^i(1-
\mu_r(A))^{n-i}$; hence $\mu_r(C) = P(\mu_r(A)) = P(Q(r))$ for all
$r \in [0,1]$.  This implies that $P \circ Q$ is a partition
polynomial. \qed

\begin{dfn}
Given partition polynomials $P,Q$, we say that $P$
\textbf{dominates} $Q$ if for some sufficiently large $n$ we can
write
\begin{eqnarray*}
P(X) &=& a_nX^n + a_{n-1}X^{n-1}(1-X) + \cdots + a_0(1-X)^n\\ Q(X)
&=& b_nX^n + b_{n-1}X^{n-1}(1-X) + \cdots + b_0(1-X)^n
\end{eqnarray*}
with $0 \leq b_i \leq a_i \leq {n \choose i}$ for each $i \leq n$.
\end{dfn}

Considering the discussion of sizes of clopen sets preceding the
above definition, we see that $P$ dominates $Q$ if and only if any
clopen set $A$ represented by $P$ has a clopen subset $B$
represented by $Q$; for if $P$ dominates $Q$ and we are given a
finite disjoint family $\U$ of cylinders with union $A$ then we
can break the members of this family into finite unions of smaller
cylinders (with a larger finite family of coordinates specified)
to obtain a family of cylinders $\U_1$, still with union $A$, that
has a subfamily the sum of whose $\mu_r$-measures is given by
$Q(r)$; the union of this subfamily is now $B$.

\section{The ideas behind the counterexample}

In this section we try to provide motivation for the
counterexample.  We will also refer back to some of the arguments
in this section in Section 3.

Suppose we are given two product measures $\mu_r$ and $\mu_s$ and
a homeomorphism $h$ such that $\mu_s = \mu_r \circ h^{-1}$.  We
observe that in this case, for each $i \in \bb{N}$ the set
$h^{-1}(A_i)$ is a clopen subset of $\O$ of $\mu_r$-measure $s$,
and, furthermore, that for any $t \in (0,s)$ the set $A_i$ has a
clopen subset of $\mu_s$-measure $t$ if and only if the set
$h^{-1}(A_i)$ has a clopen subset of $\mu_r$-measure $t$ (since
examples of the former correspond precisely to examples of the
latter under the homeomorphism $h$). Observe now that the $t$ for
which there is a subset $B$ of $A_i$ with $\mu_s(B) = t$ are
precisely those of the form $P(s)$, where $P(X)$ is a partition
polynomial that is dominated by $X$. Similarly, if the measure of
$h^{-1}(A_i)$ corresponds to the partition polynomial $Q(r)$, then
the $t$ for which there is a subset $C$ of $h^{-1}(A_i)$ with
$\mu_r(C) = t$ are precisely those of the form $R(r)$, where
$R(X)$ is a partition polynomial that is dominated by $Q(X)$.

Next we observe that a partition polynomial $P(X)$ is dominated by
$X$ if and only if it is of the form $X \cdot P_1(X)$.  Indeed, if
\[P(X) = \sum_{i=0}^na_iX^i(1-X)^{n-i}\]
and $P(X)$ is dominated by
\[X = \sum_{i=1}^n{n-1 \choose i-1}X^i(1-X)^{n-i}\]
then we must have $a_0 = 0$ and $a_i \leq {n-1 \choose i-1}$, so
we can divide $P(X)$ through by $X$ and are still left with a
partition polynomial.

We have proved the following lemma:

\begin{lem}\label{lem:simplepoly}
The fractions that can be written in the form
\[\frac{R(r)}{Q(r)} = \frac{\mu_r(C)}{\mu_r(h^{-1}(A_i))},\]
where $C$ is a clopen subset of $h^{-1}(A_i)$ (and so $R(r)$, its
$\mu_r$-measure, is given by a partition polynomial that is
dominated by $Q(X)$ evaluated at $X=r$) are precisely the
fractions of the form
\[\frac{P(s)}{s} = P_1(s)\]
where $P(X)$ is a partition polynomial dominated by $X$ and so
$P_1(X)$ is also a partition polynomial. \qed
\end{lem}

Now, since $r$ and $s$ are binomially equivalent, any value can be
written as a partition polynomial in one if and only if can be in
the other.  Thus, to find $r$ and $s$ that are binomially
equivalent but such that $\mu_r$ and $\mu_s$ are not homeomorphic,
it would suffice to find such $r$ and $s$ such that whenever $C
\subseteq \O$ has $\mu_r(C) = s$ (and so $C$ is a candidate for
any of the inverse images $h^{-1}(A_i)$) and has partition
polynomial $Q(r)$, there is a partition polynomial $R(X)$
dominated by $Q(X)$ such that the fraction $R(r) / Q(r)$ cannot be
written as a partition polynomial in $r$.  Informally, any clopen
subset $C$ of $\O$ of $\mu_r$-measure $s$ has to be so complicated
that it must contain a further clopen subset $D$ such that the
fraction $\mu_r(D) / \mu_r(C)$ is not itself the measure of any
clopen subset of $\O$.

For example, let us suppose we have the binomial reduction
\[s = 2r(1-r),\ \ \ r = F(s)\]
for some partition polynomial $F(r)$, and yet we know that no
partition polynomial in $r$ (or, equivalently, in $s$) can take
the value $\frac{1}{2}$.  (At this stage we merely speculate that
we can impose this latter condition by selecting a suitable $F$.)
Let $C$ denote the subset $\langle 1,0\rangle \cup \langle
0,1\rangle$. Then $\mu_r(C) = 2r(1-r) = s$, but no homeomorphism
$h$ sending $\mu_r$ to $\mu_s$ can be such that $C = h^{-1}(A_i)$,
because the subset $\langle 1,0 \rangle$ of $C$ has measure
$\frac{1}{2}$ that of $C$, but we know that $\frac{1}{2}$ cannot
be written as a partition polynomial.

How might we show that $\frac{1}{2}$ cannot be written as a
partition polynomial in $s$?  Were we able to write it as such, we
should have some $n \geq 1$ and some $a_0,a_1,\ldots,a_n$ with $0
\leq a_i \leq {n \choose i}$ for each $i$ such that
\[\sum_{i=0}^na_is^i(1-s)^{n-i} = \frac{1}{2}.\]
Since also
\[\sum_{i=0}^n{n \choose i}s^i(1-s)^{n-i} = (s + (1-s))^n = 1,\]
we can subtract the first of these equations from the second to
find that
\[\sum_{i=0}^na_is^i(1-s)^{n-i} = \sum_{i=0}^n\left({n \choose i} - a_i\right)s^i(1-s)^{n-i}.\]
Now we observe that in the partition polynomials on the two sides
of the above equation the coefficients of $(1-s)^n$ will either be
1 on the left and 0 on the right or 1 on the right and 0 on the
left, depending as $a_0$ is 1 or 0.  Dividing the above equation
by $r^n$ and writing $\beta = (1-r)/r$ we obtain
\[\sum_{i=0}^na_i\beta^{n-i} = \sum_{i=0}^n\left({n \choose i} - a_i\right)\beta^{n-i}.\]
Subtracting one side from the other now gives a polynomial in
$\beta$ with integer coefficients and leading term $\beta^n$.
Crucially, this requires $\beta$ to be an algebraic integer.  Thus
it will suffice to choose $F$ such that $\beta = 1/r - 1$ is
\emph{not} an algebraic integer, or, equivalently, that $1/r$ is
not an algebraic integer.  It is not hard to select such an $F$
(we will do so a little later).

Unfortunately this example does not furnish us with a
counterexample because there may be a suitable homeomorphism $h$
such that $h^{-1}$ takes each $A_i$ to a quite different clopen
subset of $\O$, one which, unlike $C$, does not have a further
clopen subset of $\frac{1}{2}$ its size.  There is no obvious way
of showing that this could not happen; in general, for $r$
algebraic, there will be many very different clopen subsets of
$\O$ of a given $\mu_r$-measure (as long as that measure is
possible at all). We do not have a way of tackling this problem;
in our actual construction we will have to take a more subtle
route instead. We will, nevertheless, still rely on ensuring that
$1/r$ is not an algebraic integer.

\section{Construction of the counterexample}

Let us choose $r$ and $s$ such that
\[s = 2r(1-r),\ \ \ \hbox{and}\ \ \ r = 3s(1-s)^2 + 3s^2(1-s).\]
It is easy to check that such $r$ and $s$ exist. Indeed, writing
$G(r) = 2r(1-r)$ and $F(s) = 3s(1-s)^2 + 3s^2(1-s)$, we see that
$F(G(\frac{1}{2})) = \frac{3}{4} > \frac{1}{2}$ and $F(G(1)) = 0 <
1$, so the intermediate value theorem assures us that suitable $r$
and $s$ exist in $(0,1)$. Many other choices are possible; we have
taken these for simplicity.  Note that $r$ and $s$ are binomially
equivalent.

\begin{lem}\label{lem:nonint}
The number $\beta = 1/r - 1$ is not an algebraic integer.
\end{lem}

\noindent\textbf{Proof}\hspace{5pt} Clearly it suffices to show that
$\alpha = 1/r$ is not an algebraic integer.  Note that $F(s)$
actually equals $3s(1-s)$ (although we have deliberately written it
as a partition polynomial above).  We have (upon substituting
$2r(1-r)$ in place of $s$):
\begin{eqnarray*}
r &=& F(s) = F(2r(1-r)) = 3(2r(1-r))(1 - 2r(1-r))\\ &=& 6r - 18r^2
+ 24r^3 - 12r^4.
\end{eqnarray*}
Re-arranging and dividing by $r$ we obtain
\[5 - 18r + 24r^2 - 12r^3 = 0,\]
and so, dividing by $r^3$,
\[5\alpha^3 - 18\alpha^2 + 24\alpha - 12 = 0.\]
By Eisenstein's criterion using the prime 3, this is irreducible;
since its leading coefficient is not $\pm 1$, $\alpha$ is not an
algebraic integer. \qed

\vspace{10pt}

We remark that for the proof of the above lemma we needed
precisely that $F$ satisfies the following conditions ($G$ is
fixed to be the one used above):
\begin{enumerate}
\item $F$ has no constant term (when written in the usual form for polynomials);
\item $F$ has a non-zero linear term (when written in the usual form for polynomials);
\item there is some prime $p \not= 2$ such that $p$ divides all
the coefficients of $F$ but $p^2$ does not.
\end{enumerate}
In the above case $p = 3$.  In general the above algebraic
manipulation yields a polynomial equation for $\alpha$ with
leading coefficient not $\pm 1$ but congruent to $1 \mod p$, with
all subsequent coefficients divisible by $p$ and with constant
term divisible by $p$ but not $p^2$, and so once again $\alpha$ is
not an algebraic integer.

\begin{theorem}
With $r$ and $s$ as above, $\mu_r$ and $\mu_s$ are not
homeomorphic.
\end{theorem}

\noindent\textbf{Proof}\hspace{5pt} Suppose, for sake of
contradiction, that there is a homeomorphism $h$ sending $\mu_r$ to
$\mu_s$.  We will construct from it an integral equation satisfied
by $\beta$, contradicting Lemma \ref{lem:nonint} above.

Consider the sets $B_i = h^{-1}(A_i)$. Since $h$ is a
homeomorphism and the sets $A_i$ and their complements generate
the topology of $\O$, the same is true of the sets $B_i$. Consider
now the points $x_j = (0,0,\ldots,0,1,0,\ldots)$ with
$j^{\rm{th}}$ coordinate equal to 1 and all others 0.  Since the
sets $B_i$ and their complements generate the topology, they must
separate these points; therefore there are $i,j$ such that $x_j
\in B_i$. Since $B_i$ is open, it follows that in fact there is
some cylinder $\langle e \rangle$ with $x_j \in \langle e \rangle
\subseteq B_i$.  By refining $\langle e \rangle$ further, we may
assume that $\langle e \rangle$ depends on the $j^{\rm{th}}$
coordinate (among others), and so must specify that this
coordinate be 1. Thus we have $\mu_r(\langle e \rangle) =
r(1-r)^{m+1}$ for some $m \geq 0$.

By the argument that proved Lemma \ref{lem:simplepoly}, we deduce
that there is some partition polynomial $K$ such that
\[\frac{\mu_r(\langle e \rangle)}{s} = \frac{r(1-r)^m}{2r(1-r)} = K(s) = K(2r(1-r)),\]
and so
\[r(1-r)^{m+1} = 2r(1-r)K(2r(1-r)).\]
Cancelling $r(1-r)$ (this is fine because $r \not= 0,1$) yields:
\[(1-r)^m = 2K(2r(1-r)).\]
Now since both $2X(1-X)$ and $K(X)$ are partition polynomials, so
is their composition $K(2X(1-X))$, by Lemma \ref{lem:polycomp};
hence, in particular, we may write
\[K(2r(1-r)) = c_kr^k + c_{k-1}r^{k-1}(1-r) + \cdots + c_0(1-r)^k\]
with $c_0$ equal to 0 or 1.  Our equation for $r$ becomes
\[(1-r)^m = 2(c_kr^k + c_{k-1}r^{k-1}(1-r) + \cdots + c_0(1-r)^k).\]
If $k \leq m$ we can repeatedly replace terms $T$ on the right
hand side with sums of terms $Tr + T(1-r)$ and so assume that $k =
m$; if, on the other hand, $k > m$ then we let $p = k - m$ and
write instead
\begin{eqnarray*}&&(1-r)^k + pr(1-r)^{k-1} + \cdots +
r^p(1-r)^m\\ &&\hspace{10pt} = 2(c_kr^k + c_{k-1}r^{k-1}(1-r) +
\cdots + c_0(1-r)^k).\end{eqnarray*} Dividing this equation by
$r^{m+p}$ gives
\[\beta^k + p\beta^{k-1} + \cdots + \beta^m = 2(c_k + c_{k-1}\beta + \cdots + c_0\beta^k).\]
Here the highest term is in $\beta^k$, with coefficient 1 on the
left hand side and either 0 or 2 on the right hand side; either way
we obtain an integral equation for $\beta$, and so the desired
contradiction. \qed

\vspace{10pt}

Following on from the remark after Lemma \ref{lem:nonint}, we
observe that the precise form of $F$ did not enter the above proof
at all.  Thus the counterexample given here is in no way special;
others can be constructed from this $G$ and any $F$ satisfying the
conditions listed after that Lemma (subject only to the further
requirement that roots $r$ and $s$ exist in $(0,1)$ at all; often
this is clear from the intermediate value theorem, but if not it
can be proved using Brouwer's fixed point theorem, as in
\cite{Mauld}).

\section*{Acknowledgements} The above work was carried out under
a summer research studentship funded by Trinity College, Cambridge
over the Long Vacation period 2005.  My thanks go to Dr I. Leader
for his advice and the official supervision of the project.

%\bibliographystyle{alpha}
%\bibliography{bibfile}

\noindent \textsc{Department of Mathematics, University of California at Los Angeles, Los Angeles CA 90095, USA}\\
{\em email: }\verb|timaustin@math.ucla.edu|

\end{document}